\documentclass{amsart}

\usepackage[utf8]{inputenc}
\usepackage[T1]{fontenc}
\usepackage{lmodern}
\usepackage{scalerel}
\usepackage[only,boxbslash]{stmaryrd}
\usepackage[colorlinks=true,citecolor=magenta,urlcolor=blue]{hyperref}

\newtheorem{theorem}{Theorem}
\newtheorem{lemma}{Lemma}
\newtheorem{proposition}{Proposition}

\DeclareMathOperator{\Tr}{tr}
\newcommand{\dif}{\mathrm d}
\newcommand{\nDelta}{{\mathord{\scalerel*{\boxbslash}{gX}}}}

\begin{document}

\title{Kinetic Dyson Brownian motion}
\author{Pierre Perruchaud}
\address{Department of Mathematics, University of Notre Dame, Notre Dame, IN 46556, USA.}
\email{pperruch@nd.edu}

\begin{abstract}
We study the spectrum of the kinetic Brownian motion in the space of $d\times d$ Hermitian matrices, $d\geq2$. We show that the eigenvalues stay distinct for all times, and that the process $\Lambda$ of eigenvalues is a kinetic diffusion (i.e. the pair $(\Lambda,\dot\Lambda)$ of $\Lambda$ and its derivative is Markovian) if and only if $d=2$. In the large scale and large time limit, we show that $\Lambda$ converges to the usual (Markovian) Dyson Brownian motion under suitable normalisation, regardless of the dimension.
\end{abstract}

\maketitle

\section{Introduction}

In the space $\mathcal H_d$ of $d\times d$ complex Hermitian matrices, there is a natural Brownian process $W$, whose covariance is given by the Hilbert--Schmidt norm. It turns out, as was first observed by Dyson in 1962 \cite{Dyson}, that the spectrum of $W$ is Markovian: the law of the spectrum of $s\mapsto W_{t+s}$ depends on $W_t$ only through its spectrum. In this paper, we show that for a natural smoothing of the Brownian motion, the so-called Kinetic Brownian motion, the situation can be different.

Here, a kinetic motion with values in $\mathcal H_d$ is a process of regularity $\mathcal C^1$, such that the couple $(H,\dot H)$ of the position and associated velocity is Markovian. Kinetic Brownian motion is the kinetic motion whose velocity $\dot H$ is a standard Brownian motion on the unit sphere. We say that it is a smoothing of Brownian motion because in the large scale limit, $H$ looks like a Brownian motion; namely, the law of the process $t\mapsto\frac1LH_{L^2t}$ converges to that of $W$ as $L\to\infty$, up to a factor $4/d^2(d^2-1)$. In other words, $\frac1LH_{L^2t}$ is very similar to a Brownian motion, with the important difference that it is actually $\mathcal C^1$. For a proof of this convergence, see \cite[Theorem 1.1]{LiKBM} or \cite[Proposition 2.5]{ABT}. See references below for more on kinetic Brownian motion.

Define $\Lambda_t$ as the vector of eigenvalues of $H_t$, the order being irrelevant provided it depends continuously on time. Then a classical application of the inverse function theorem shows that $\Lambda$ has to be $\mathcal C^1$ whenever the eigenvalues of $H_t$ are distinct; in particular, $\Lambda$ cannot be Markovian, otherwise it would be deterministic. The next natural hope would be for the process $(\Lambda,\dot\Lambda)$ to be Markovian. The main objective of this paper is to prove it is the case if and only if $d=2$.

\begin{theorem}
\label{thm:main}
Let $(H,\dot H)$ be a kinetic Brownian motion on the space $\mathcal H_d$ of $d\times d$ complex Hermitian matrices ($d\geq2$), and $0\leq\tau\leq\infty$ the first time $H$ has multiple eigenvalues. Let $\Lambda$ be the process of eigenvalues of $H$, seen as a continuous process from an interval of $\mathbb R_+$ to $\mathbb R^d$.

Then $\tau = \infty$ whenever $H_0$ has distinct eigenvalues. Moreover, $(\Lambda_t,\dot\Lambda_t)_{0\leq t<\tau}$ is well-defined, and is Markovian if and only if $d=2$.
\end{theorem}

In the case $d=2$, the stochastic differential equations describing the eigenvalues of $H$ are given in Lemma \ref{lem:casesond}, Section \ref{ssec:criterion}. In the general case, we introduce in Section \ref{ssec:LambdaA} a subdiffusion $(\Lambda,A)$ of $(H,\dot H)$. Since this process is Markovian in any dimension (see Lemma \ref{lem:MarkovLambdaA} for defining equations), one may want to consider it as a suitable approach to kinetic Dyson Brownian motion, rather than the more restrictive $(\Lambda,\dot\Lambda)$.

As discussed above, it is known that the process $t\mapsto\frac1LH_{L^2t}$ converges in law to a Brownian motion. From this known fact, we will deduce that $\Lambda$ converges to a standard Dyson Brownian motion, in the following sense.

\begin{proposition}
\label{prop:homogenisation}
Let $(H,\dot H)$ be a kinetic Brownian motion on the space $\mathcal H_d$ of $d\times d$ complex Hermitian matrices ($d\geq2$) such that $H_0$ has distinct eigenvalues almost surely. Let $\Lambda$ be the process of eigenvalues of $H$ in non-decreasing order, seen as a continuous process from an interval of $\mathbb R_+$ to $\mathbb R^d$. Let $D$ be the process of eigenvalues of a standard Brownian motion in $\mathcal H_d$ starting at zero, with the same conventions.

Then we have the following convergence in law as $L$ goes to infinity:
\[  (t\mapsto \Lambda_{L^2t})_{0\leq t\leq1}\xrightarrow{\mathcal L}\frac4{d^2(d^2-1)}\cdot D.  \]
\end{proposition}

We give a sketch of proof of Theorem \ref{thm:main} in Sections \ref{ssec:tau} to \ref{ssec:criterion}, using various lemmas proved in Part \ref{sec:proofs}. Proposition \ref{prop:homogenisation} is proved in Section \ref{ssec:homogenisation}.

Kinetic Brownian motion was first introduced by Li in \cite{LiIntertwining}, where Theorem 4.3 proves a stronger convergence theorem to Brownian motion than stated above. A self-contained proof by the same author appeared later in \cite{LiKBM}. This result was generalised by Angst, Bailleul and Tardif \cite{ABT} and the author \cite{Perruchaud}. Kinetic Brownian motion has been given different names in the literature, for instance velocity spherical Brownian motion by Baudoin and Tardif \cite{BaudoinTardif} or circular Langevin diffusion by Franchi \cite{Franchi} in the context of heat kernels. See also \cite{ABP} for considerations in an infinite dimensional setting.

The question of existence and behaviour of kinetic Dyson Brownian motion was raised by Thierry Lévy during the author's PhD defence. It turned out to be a beautiful endeavour, and the latter would like to thank the former for this suggestion.

\section{Definitions and Proof Outline}

Let $\mathcal H_d$ be the space of $d\times d$ complex Hermitian matrices. We assume $d\geq2$, and endow it with the Hilbert--Schmidt inner product:
\[ \|H\|_\mathcal{H}^2:=\Tr(H^2)=\Tr(H^*H)=\sum_{ij}|H_{ij}|^2, \]
where $H^*$ is the conjugate transpose of $H$, $(H_{ij})_{ij}$ are the coefficients of $H$, and $\langle\cdot,\cdot\rangle$ is the Hermitian product on $\mathbb C^d$, with the convention that it is linear in its second argument. It is isometric to the standard Euclidean matrix space $\mathbb R^{d\times d}$, via
\[ M=(m_{ij}) \mapsto
   \begin{pmatrix}
     m_{11} & \frac{m_{12}+\mathsf i\,m_{21}}{\sqrt 2} & & \cdots & \frac{m_{1d}+\mathsf i\,m_{d1}}{\sqrt 2} \\
     \frac{m_{12}-\mathsf i\,m_{21}}{\sqrt 2} & m_{22} & \ddots & & \\[15pt]
     \vdots & \ddots & \ddots & \ddots & \vdots  \\[15pt]
      & & \ddots & m_{d-1,d-1} & \frac{m_{d-1,d}+\mathsf i\,m_{d,d-1}}{\sqrt 2} \\
     \frac{m_{1d}-\mathsf i\,m_{d1}}{\sqrt 2} & \cdots & & \frac{m_{d-1,d}-\mathsf i\,m_{d,d-1}}{\sqrt 2} & m_{dd}
     \end{pmatrix}. \]
A Brownian motion $W$ in $\mathcal H_d$ associated to this Euclidean structure can be described as a matrix as above, where the $m_{ij}$'s are independent (real standard) Brownian motions.

Let $\dot H$ be a standard Brownian motion on the unit sphere $\mathbb S(\mathcal H_d)$ of $\mathcal H_d$, and $H$ its integral. For instance, one may define $(H,\dot H)$ as the solution of the stochastic differential equation
\begin{align}
\label{eq:defH}
\dif H_t     & = \dot H_t\dif t, \\
\label{eq:defdotH}
\dif\dot H_t & = \underbrace{{}\circ\dif W_t - \dot H_t\langle\dot H_t,{}\circ\dif W_t\rangle_\mathcal{H}}_{\text{projection of }\circ\dif W_t\text{ on }\dot H_t^\perp}
               = \dif W_t - \dot H_t\Tr(\dot H_t^*\dif W_t) - \frac{d^2-1}2\dot H_t\dif t,
\end{align}
where $\circ\dif W$ (resp. $\dif W$) denotes the Stratonovich (resp. Ito) integral. It is defined for all times, since $\dot H$ is the solution of a SDE with smooth coefficients on a compact manifold, and $H$ is the integral of a process that is uniformly bounded. Let $0\leq\tau\leq\infty$ be the first time $H$ has multiple eigenvalues, with $\tau=\infty$ if its eigenvalues stay distinct for all times. It is a stopping time since it is the hitting time of a closed set.

In the following, we will work with diagonal matrices and matrices whose diagonal is zero; let $\mathcal H_d=\mathcal H_d^\Delta\oplus\mathcal H_d^\nDelta$ be the associated decomposition. Note that it is actually an orthonormal decomposition. Let us also write $\mathfrak u_d$ for the space of skew-Hermitian matrices ($\mathfrak u_d$ is the Lie algebra of the group $U_d(\mathbb C)$ of unitary matrices).

The remainder of this section is a complete proof outline of Theorem \ref{thm:main}, using a few lemmas proved in the next section.

\subsection{Explosion time}
\label{ssec:tau}

By ``$\tau=\infty$ whenever $H_0$ has distinct eigenvalues'', we mean that the event
\[ \{H_0\text{ has distinct eigenvalues and }\tau<\infty\} \]
has measure zero. It is known, see references in Section \ref{ssec:discriminant}, that the subset of $\mathcal H_d$ consisting of matrices with multiple eigenvalues is contained in a finite collection of submanifolds of codimension 3. Then it is enough to prove the following result, of independent interest.

\begin{proposition}[proved in Section \ref{ssec:codim2}]
\label{prop:codim2}
Let $M$ be a complete Riemannian manifold, and $N\subset M$ a submanifold of codimension at least 2. Let $(H,\dot H)$ be a kinetic Brownian motion in $M$, defined for all times $t\geq0$. Then the event
\[ \{ H_t\in N\text{ for some }t>0\} \]
has measure zero.
\end{proposition}

\subsection{Diagonalisation}

Because $H_t$ is Hermitian, there exists for all $t$ a unitary matrix $U_t$ such that $U_t^*H_tU_t$ is diagonal. Abstract geometric arguments show that it is possible to show that for a fixed realisation of $H$, we can find a $U$ with regularity $\mathcal C^1$, at least as long as the eigenvalues stay distinct. However, we would like $U$ to be described by an explicit stochastic differential equation. Let us look for a candidate.

Given a $\mathcal C^1$ process $U$ with unitary values, we call
\[ \dot u_t := U_t^{-1}\dot U_t = U_t^*\dot U_t \]
its derivative, seen in the Lie algebra of $U_d(\mathbb C)$: $\dot u_t\in\mathfrak u_d$. If we define the $\mathcal C^1$ process
\begin{equation}
\label{eq:defLambda}
\Lambda:t\mapsto U_t^*H_tU_t,
\end{equation}
then $H_t$ will be diagonal in the frame $U_t$ if and only if $\Lambda_t\in\mathcal H^\Delta$. It means that we are looking for a $U$ such that the derivative $\dot\Lambda_t$ stays in $\mathcal H^\Delta$. We have
\[ \dot\Lambda_t
 = \big(\dot U_t^*H_tU_t + U_t^*H_t\dot U_t\big) + U_t^*\dot H_tU_t
 = (\dot u_t^*\Lambda_t + \Lambda_t\dot u_t) + U_t^*\dot H_tU_t. \]
Assuming $\Lambda$ is indeed diagonal, and since $\dot u_t$ is skew-Hermitian, the coefficients of the first term are
\begin{equation}
\label{eq:productdotuLambda}
   (\dot u_t^*\Lambda_t + \Lambda_t\dot u_t)_{ij}
 = (\dot u_t^*)_{ij}(\Lambda_t)_{jj} + (\Lambda_t)_{ii}(\dot u_t)_{ij}
 = \big((\Lambda_t)_{ii} - (\Lambda_t)_{jj}\big)(\dot u_t)_{ij}.
\end{equation}
In other words, for $\Lambda$ to stay diagonal, we have no choice for the off-diagonal coefficients of the velocity $\dot u_t$: setting
\begin{equation}
\label{eq:defA}
A:t\mapsto U_t^*\dot H_tU_t,
\end{equation}
they have to be
\[  -\frac{(U_t^*\dot H_tU_t)_{ij}}{(U_t^*H_tU_t)_{ii}-(U_t^*H_tU_t)_{jj}}
  = -\frac{(A_t)_{ij}}{(\Lambda_t)_{ii}-(\Lambda_t)_{jj}}. \]
It turns out that this choice works, as we will see.

For the sake of conciseness, define $\dot u(\Lambda,A)$ as
\begin{equation}
\label{eq:defdotu}
\big(\dot u(\Lambda,A)\big)_{ij} :=
\begin{cases}
  \displaystyle-\frac{A_{ij}}{\Lambda_{ii}-\Lambda_{jj}} & \text{ for } i\neq j, \\
  0 & \text{ else}
\end{cases}
\end{equation}
whenever $\Lambda$ has distinct diagonal entries. As long as $\Lambda$ stays diagonal with distinct eigenvalues, $u(\Lambda,A)$ stays well-defined.

\begin{lemma}[proved in Section \ref{ssec:proofL12}]
\label{lem:defULambda}
Let $(H_t,\dot H_t)_{t\geq0}$ be a kinetic Brownian motion on $\mathcal H_d$, and $0\leq\tau\leq\infty$ the first time $H$ has multiple eigenvalues. Let $U_0$ be a (random) unitary matrix such that $U_0^*H_0U_0$ is diagonal, and define $U_t$ as the solution of
\begin{align}
\label{eq:defU}
\dif U_t &= U_t\dot u_t\dif t, &
\dot u_t &= \dot u\big(U_t^*H_tU_t,U_t^*\dot H_tU_t\big),
\end{align}
where $\dot u(\Lambda,A)$ is defined in equation \eqref{eq:defdotu}.

Then $U_t$ is defined for all $0\leq t<\tau$, and $U_t^*H_tU_t$ is diagonal for all such $t$.
\end{lemma}

In particular, it means that $\Lambda$ is the process of eigenvalues of $H$, so the potential kinetic Dyson Brownian motion is $(\Lambda,\dot\Lambda)$. If one wishes, we can take $U_0$ so that the diagonal entries of $U_0^*H_0U_0$ are in a given order, for instance non-decreasing. Then, using the fact that the eigenvalues of $H_t$ stay distinct for all $t>0$ (see Section \ref{ssec:tau}), we see that the diagonal entries $U_t^*H_tU_t$ stay in the same order. In the following, we will not need nor assume that $U_0$ satisfies this property.

\subsection{The Markovian process \texorpdfstring{$(\Lambda,A)$}{(Λ,A)}}
\label{ssec:LambdaA}

From $(H,\dot H)$, we constructed $U$ using equation \eqref{eq:defU}. We define then $\Lambda$ and $A$ as in equations \eqref{eq:defLambda} and \eqref{eq:defA}. By Lemma \ref{lem:defULambda} we know that $\Lambda$ is in fact diagonal. We also notice that $A$ takes values in the sphere $\mathbb S(\mathcal H_d)$, since $\dot H$ does, and conjugation by a fixed unitary matrix is an isometry.

Then we see that the triple $(U_t,\Lambda_t,A_t)$ satisfies the system of equations
\begin{align*}
\dif U_t & = U_t\dot u_t\dif t \\
\dif \Lambda_t & = (\dot u_t^*\Lambda_t + \Lambda_t\dot u_t)\dif t + A_t\dif t \\
\dif A_t & = (\dot u_t^*A_t + A_t\dot u_t)\dif t + U_t^*\dif\dot H_tU_t
\end{align*}
for $0\leq t<\tau$, where $\dot u_t = \dot u(\Lambda_t,A_t)$. These $\dot u(A,\Lambda)$ and $\dot u_t$ are still the same, defined respectively in equations \eqref{eq:defdotu} and \eqref{eq:defU}; we are merely emphasising the fact that $\dot u_t$ depends on $(H,\dot H)$ only through $(\Lambda,A)$. Note that the above is really a stochastic differential equation describing $(U,\Lambda,A)$ with a driving noise $\dot H$, rather than an abstract functional of the couple $(H,\dot H)$. If one wishes, one can see $(U,\Lambda,A)$ as a continuous process defined for all times, over the one-point compactification $\mathcal U\cup\{\dagger\}$ of the open set $\mathcal U$ of triples $(U,\Lambda,A)\in U_d(\mathbb C)\times\mathcal H_d^\Delta\times\mathbb S(\mathcal H_d)$ such that $\Lambda$ has distinct eigenvalues, with the convention that $\dagger$ is absorbing. However, since $\tau$ can only be zero or infinity, the benefit of doing so is rather small.

We are interested in the dynamics of $(\Lambda,\dot\Lambda)$. In fact, this pair is nothing but $(\Lambda,\pi^\Delta(A))$, where the second term is the projection of $A$ on $\mathcal H_d^\Delta$, roughly speaking its diagonal. Indeed, this is a direct consequence of the fact that in the driving equation for $\dif\Lambda_t$, the first term is zero on the diagonal as seen in \eqref{eq:productdotuLambda}. As explained in the following lemma, it turns out that $(\Lambda,A)$ is Markovian.

\begin{lemma}[proved in Section \ref{ssec:proofL12}]
\label{lem:MarkovLambdaA}
Let $(H,\dot H)$ be a kinetic Brownian motion in $\mathcal H_d$, $0\leq\tau\leq\infty$ the first time $H$ has multiple eigenvalues, and define $U$, $\Lambda$ and $A$ as in equations \eqref{eq:defU}, \eqref{eq:defLambda} and \eqref{eq:defA}. Then, up to enlarging the underlying probability space, there exists a standard Brownian motion $B$ with values in $\mathcal H_d$ such that
\begin{align*}
\dif \Lambda_t
& = \pi^\Delta(A_t)\dif t \\
\dif A_t
& = \big(\dot u_t^*A_t + A_t\dot u_t\big)\dif t
    + \dif B_t - A_t\Tr\big(A_t^*\dif B_t\big)
    - \frac{d^2-1}2A_t\dif t
\end{align*}
for all $0\leq t<\tau$, where $\dot u_t = u_t(\Lambda_t,A_t)$ as defined in equation \eqref{eq:defdotu} and $\pi^\Delta
$ is the projection on the space of diagonal matrices. In particular, the process $(\Lambda,A)$ is a Markovian process.
\end{lemma}
Here, we could replace $\dot u_t$ everywhere by its actual expression, which makes it clear that $(\Lambda,A)$ satisfies a self-contained stochastic differential equation. We can again see $(\Lambda,A)$ as a process defined for all times, over the one-point compactification of an appropriate open set of $\mathcal H_d^\Delta\times\mathbb S(\mathcal H_d)$.

Since $(\Lambda,A)$ is Markovian in all dimensions whereas (as we will show later) $(\Lambda,\dot\Lambda)$ is not, we made the remark in the introduction that one might view the former as a natural definition of kinetic Dyson Brownian motion. Instead of containing only the information about the derivative of the eigenvalues of $\Lambda$ as $\dot\Lambda$ might do, $A$ is the matrix $\dot H$ as seen in a referential that makes $H$ diagonal; in particular, it contains some hints about the motion of the eigenspaces in relation to each other. As we will see below, at least some of this additional information is needed to describe the motion of $\Lambda$ entirely, at least in dimension $d\geq3$.

\subsection{A criterion for a Markovian kinetic Dyson Brownian motion}
\label{ssec:criterion}

It is worth noticing that the equation for $A$ describes a Brownian motion with drift on a sphere; see for instance the definition of $\dot H_t$ in \eqref{eq:defdotH}, describing a Brownian motion without drift. In fact, if it were not for the first term $(\dot u^*A+A\dot u)\dif t$, $A$ would be precisely a standard Brownian motion on the unit sphere of $\mathcal H_d$.

But it is known, and not too difficult to see, that the projection of a spherical Brownian motion $X$ is Markovian. Indeed, once one fixes a subset of coordinates $(X^0,\ldots,X^k)$ of norm $r$, then the remaining coordinates $(X^{k+1},\ldots,X^n)$ can always be reduced to $\big(\sqrt{1-r^2},0,\ldots,0\big)$, up to a rotation fixing the first coordinates; see Section \ref{ssec:sphericalBM} for details and references. In particular, if we continue to ignore this drift term, it would be clear at this point that $(\Lambda,\dot\Lambda)=(\Lambda,\pi^\Delta(A))$ is Markovian. So any obstruction for $(\Lambda,\dot\Lambda)$ to be Markovian must come from the additional term $(\dot u^*A+A\dot u)$. The following lemma describes the situation with a precise criterion.

\begin{lemma}[proved in Section \ref{ssec:proofL3}]
\label{lem:MarkovPhi}
Define $\Phi=\Phi^{(d)}:\mathcal H_d^\Delta\times\mathbb S(\mathcal H_d)\to \mathcal H_d^\Delta$ by
\[  \Phi(\Lambda,A)
 := \pi^\Delta(\dot u^*A+A\dot u), \]
where $\pi^\Delta$ is the projection on $\mathcal H_d^\Delta$ and $\dot u = \dot u(\Lambda,A)$ as defined in \eqref{eq:defdotu}.

Let $(H,\dot H)$ be a kinetic Brownian motion in $\mathcal H_d$, and define $\Lambda$ the continuous process of its eigenvalues. Then $(\Lambda,\dot\Lambda)$ is Markovian if and only if $\Phi$ factors through $(\Lambda,A)\mapsto(\Lambda,\pi^\Delta(A))$.
\end{lemma}
There is a concise expression for the coefficients of $\Phi(\Lambda,A)$. It is of course zero out of the diagonal, and we have
\begin{equation}
\label{eq:concisePhi}
\Phi(\Lambda,A)_{ii}
 = \sum_{j\neq i}\Big(- \frac{\overline A_{ji}}{\Lambda_{jj}-\Lambda_{ii}}\cdot A_{ji}
                      - A_{ij}\cdot\frac{A_{ji}}{\Lambda_{jj}-\Lambda_{ii}}\Big)
 = 2\sum_{j\neq i}\frac{|A_{ij}|^2}{\Lambda_{ii}-\Lambda_{jj}}.
\end{equation}
In dimension 2, since $|A| = 1$, we get directly
\[ \Phi(\Lambda,A)_{11}
 = - \Phi(\Lambda,A)_{22}
 = \frac{|A_{12}|^2+|A_{21}|^2}{\Lambda_{11}-\Lambda_{22}}
 = \frac{1-|A_{11}|^2-|A_{22}|^2}{\Lambda_{11}-\Lambda_{22}}. \]
This depends only on $\Lambda$ and the on-diagonal coefficients of $A$, so $(\Lambda,\dot\Lambda)$ is indeed Markovian. In dimension $d\geq3$, it is not obvious that one could use a similar trick, and in fact we can show that it is not possible and that the process is not Markovian.

In Sections \ref{ssec:dim2} and \ref{ssec:dim3+} we carry out the computations in dimension $d=2$ and $d\geq3$ respectively, and we conclude as follows.

\begin{lemma}[proved in Sections \ref{ssec:dim2} and \ref{ssec:dim3+}]
\label{lem:casesond}
~

\begin{itemize}
\item If $d=2$, the eigenvalues $\lambda,\mu$ of $H$ make up a kinetic diffusion, and satisfy the equations
\begin{align*}
\dif\lambda_t &= \dot\lambda_t\dif t, &
\dif\dot\lambda_t &= +\frac{1-\dot\lambda_t^2-\dot\mu_t^2}{\lambda_t-\mu_t}\dif t
                     + \dif M^\lambda_t
                     - \frac{d^2-1}2\dot\lambda_t\dif t, \\
\dif\mu_t &= \dot\mu_t\dif t, &
\dif\dot\mu_t &= -\frac{1-\dot\lambda_t^2-\dot\mu_t^2}{\lambda_t-\mu_t}\dif t
                     + \dif M^\mu_t
                     - \frac{d^2-1}2\dot\mu_t\dif t,
\end{align*}
where $M^\lambda$ and $M^\mu$ are martingales with brackets
\begin{align*}
\dif\langle M^\lambda,M^\lambda\rangle_t &= \big(1 - \dot\lambda_t^2\big)\dif t, &
\dif\langle M^\mu,M^\mu\rangle_t &= \big(1 - \dot\mu_t^2\big)\dif t,
\end{align*}
\[ \dif\langle M^\lambda,M^\mu\rangle_t = - \dot\lambda_t\dot\mu_t\dif t. \]

\item If $d\geq3$, $\Phi^{(d)}$ does not factor, and the process $(\Lambda,\dot\Lambda)$ is not Markovian.
\end{itemize}
\end{lemma}

It is tempting to try to salvage some weaker Markovian properties from $(\Lambda,\dot\Lambda)$ when $d\geq3$. In Section \ref{ssec:dim3+}, we show that in some sense, it is not locally Markovian anywhere.

\subsection{Homogenisation}
\label{ssec:homogenisation}

As stated in the introduction, if we write $H^L$ for the normalised process $t\mapsto\frac1LH_{L^2t}$, then $(H^L_t)_{0\leq t\leq1}$ converges in law to a standard Brownian motion, up to a constant scaling factor $4/d^2(d^2-1)$; see \cite[Theorem 1.1]{LiKBM} or \cite[Proposition 2.5]{ABT}. In this section we prove Proposition \ref{prop:homogenisation}, namely that the process $\Lambda$, although not Markovian, is somehow Markovian in large scales, in the sense that a similar limit converges to a Dyson Brownian motion.

The map $\mathcal H_d\to\mathcal H_d^\Delta$ sending a matrix $H$ to the matrix $\Lambda$ whose diagonal entries are the eigenvalues of $H$ with multiplicities according to a chosen order (e.g. non-decreasing) is continuous. One can see this as follows. First, the spectral measure $\frac1d\sum\delta_\lambda$ is a continuous function of $M\in M_d(\mathbb C)$, where $\delta_\lambda$ is the Dirac mass at $\lambda$, and $\lambda$ ranges over the spectrum of $M$. Then, since $H$ has real eigenvalues, we only need continuity of the map sending a spectral measure on $\mathbb R$ to the ordered vector of its atoms, which follows directly from the easy fact that the smallest atom depends continuously on the measure.


This means that $\Lambda_t$ can more or less be described as a continuous function of $H_t$. We say more or less, because the ordering of the eigenvalues depends on our choice of $U_0$: if it was chosen such that the eigenvalues of $U_0^*H_0U_0$ are in non-decreasing order, which we can always impose, then the values of $\Lambda_t$ along the diagonal would stay in non-decreasing order (since the eigenvalues of $H_t$ are distinct for all $t>0$, see Section \ref{ssec:tau}). Suppose we chose $U_0$ as such, and define $\Lambda:\mathcal H_d\to\mathcal H_d^\Delta$ as the map sending a given matrix $H$ to the diagonal matrix whose entries are the eigenvalues of $H$ in non-decreasing order; it should lead to no confusion of notation, since we then have $\Lambda_t=\Lambda(H_t)$. In particular, if one defines $\Lambda^L:t\mapsto\frac1L\Lambda_{L^2t}$, then $\Lambda^L=\Lambda\circ H^L$. Such operations preserve convergence in law, so given a standard Brownian motion $W$ in $\mathcal H_d$, we have the following convergence in law:
\[  (H^L_t)_{0\leq t\leq1}\xrightarrow{\mathcal L}\frac4{d^2(d^2-1)}\cdot W,
    \qquad
    (\Lambda^L_t)_{0\leq t\leq1}\xrightarrow{\mathcal L}\frac4{d^2(d^2-1)}\cdot\Lambda(W).  \]
Since $\Lambda(W)$ is the spectrum of a Brownian motion in $\mathcal H_d$ in the form of a diagonal matrix, it is nothing but a Dyson Brownian motion. As stated above, $\Lambda$ looks very much like a Dyson Brownian motion at large scales. Recall that Dyson Brownian motion is Markovian, so the hidden information preventing $\Lambda$ to be Markovian (the off-diagonal coefficients of $A$, but also the derivative $\dot\Lambda$) vanishes in the limit.

It might be interesting to see if one could prove the convergence of $\Lambda^L$ towards the rescaled $\Lambda(W)$ using only the dynamics of $(\Lambda,A)$ as given in Lemma \ref{lem:MarkovLambdaA}. Although the author does not pretend it is impossible, it seems that the non-linearity in the vector field $\Phi(\Lambda,A)$ makes it more difficult to approach than the convergence of $H^L$, using for instance the methods of \cite{ABT}.

\section{Proof of the Lemmas}
\label{sec:proofs}

\subsection{Matrices with multiple eigenvalues}
\label{ssec:discriminant}

We claimed earlier that the set of Hermitian matrices with multiple eigenvalues is covered by finitely many submanifolds of codimension 3. This is all we will need for our purposes, although we can actually show an additional structure result. Since this set is a (real) algebraic set (it is the zero locus of the discriminant of the characteristic polynomial), it actually means that it is the disjoint union of at most $d^2-2$ (real analytic) submanifolds of dimensions $0$ to $d^2-3$; see for instance propositions 3.3.11 and 3.3.14 of \cite{BochnakCosteRoy}.

We sketch a proof of the covering result; the approach is carried out in \cite{Arnold} with a detailed discussion of the underlying combinatorial structure. See also \cite{RepeatedEigenvalues} and references therein. Let $d_1+\cdots+d_k=d$ be a partition of $d$ into $k$ positive integers. Consider the space $N$ of all Hermitian matrices $H$ with eigenvalues $\lambda_1\leq\cdots\leq\lambda_d$ such that the first $d_1$ eigenvalues are equal but less than the next, the following $d_2$ are equal but less than the $(d_1+d_2+1)$st, etc. For instance, if $(d_1,d_2,d_3)=(1,2,2)$, we are considering $5\times5$ matrices whose eigenvalues satisfy
\[ \lambda_1<\lambda_2=\lambda_3<\lambda_4=\lambda_5. \]
There is a one-to-one correspondence between such matrices and a choice of $k$ (orthogonal) eigenspaces and associated (real) eigenvalues. There are $k$ degrees of freedom for the choice of the eigenvalues. Using for instance the Gram--Schmidt algorithm, the choice of the eigenspaces is equivalent to the data of a flag $E_1\subset\cdots\subset E_k=\mathbb C^d$ of subspaces of respective dimensions $i_1\leq\cdots\leq i_k$. The space of these flags is known to be a manifold of complex dimension
\[ d_1(d-i_1) + d_2(d-i_2) + \cdots + d_{k-1}(d-i_{k-1})
 = \frac12\Big(d^2 - \sum_\ell d_\ell^2\Big). \]
All in all, the set of matrices satisfying this constraint is a manifold of real dimension
\[ d^2-(d_1^2-1)-\cdots-(d_k^2-1) \]
(the restriction to a space of dimension $d_\ell$ could have been any matrix of $\mathcal H_{d_\ell}$, but is instead scalar), so it has codimension at least $3$ when a given $d_\ell$ is not one. Considering all partitions of $d$ except the trivial $1+\cdots+1$, we see that the set of matrices with multiple eigenvalues is included in a finite collection of manifolds of codimension at least 3.

\subsection{Proof of Proposition \ref{prop:codim2}}
\label{ssec:codim2}

Let us turn to the proof of Proposition \ref{prop:codim2}. Let $(M,g)$ be a complete Riemannian manifold of dimension $n$, and $N\subset M$ a submanifold of codimension at least 2. We suppose $N$ is an embedded manifold without boundary, although it will be clear that the proof may be adapted to the more general case of immersed manifolds with boundary. Given a kinetic Brownian motion $(H,\dot H)$ on $M$, we want to show that the event
\[ \{ H_t\in N\text{ for some }t>0\} \]
has probability zero.

We call embedded (closed) disc of dimension $k$ a subset $D$ of $M$ for which there exists an open set $\mathcal U\supset D$ and a diffeomorphism $\phi:\mathcal U\to B_0(1)$ to the unit ball of $\mathbb R^n$ such that $\phi(D)$ is the intersection $(\mathbb R^k\times\{0\}^{n-k})\cap\overline B_0(1/2)$. Then, because $N$ is second countable, it can be covered by countably many embedded discs of codimension 2, say $N\subset\bigcup_{i\geq0}D_i$. It means that
\[   \mathbb P(H_t\in N\text{ for some }t>0)
\leq \sum_{N,i\geq0}\mathbb E\big[\mathbb P(H_t\in D_i\text{ for some }t\in[2^{-N},2^N])\big|(H_0,\dot H_0)\big]. \]
We are left to show that for a given compact interval $[a,b]\subset(0,\infty)$, starting point $(H_0,\dot H_0)$ and embedded disc $D$ of codimension 2, we have
\[ \mathbb P\left( H_t\in D\text{ for some } t\in [a,b] \right) = 0. \]

Fix some $\delta>0$, and write $D_\delta$ for the set of points at distance at most $\delta$ from $D$. Since the position process of the kinetic Brownian motion has velocity one, if we have $H_t\in D$ for a given $t>0$, then there must exist $t'\in2\delta\mathbb N$ such that $H_{t'}\in D_\delta$. It means that
\begin{equation}
\label{eq:discretetime}
      \mathbb P\left( H_t\in D\text{ for some } t\in I \right)
 \leq \sum_{\ell=\lfloor a/2\delta\rfloor}^{\lceil b/2\delta\rceil}\mathbb P(H_{2\delta\ell}\in D_\delta).
\end{equation}
We will prove that
\[ \limsup_{\delta\to0}\sup_{a/2\leq t\leq2b}\frac{\mathbb P(H_t \in D_\delta)}{\delta^2} < \infty, \]
which will show that the sum in \eqref{eq:discretetime} is bounded by a constant multiple of $\delta$, so the left hand side is zero upon taking the limit $\delta\to0$.

Let $\phi:\mathcal U\to B_0(1)$ be a map compatible with $D$, in the sense described above. It induces a diffeomorphism $T^1\phi$ from $T^1\mathcal U$ to $T^1B_0(1)\simeq B_0(1)\times\mathbb S^{n-1}$, sending $(h,\dot h)$ to
\[ \left(\phi(h),\frac{d\phi_h(\dot h)}{|d\phi_h(\dot h)|}\right). \]
Since $D$ is compact, there is some small $\varepsilon>0$ such that $D_\varepsilon\subset\mathcal U$. The fact that $(H,\dot H)$ is a hypoelliptic diffusion means that the density of $(H_t,\dot H_t)$ is smooth for any given $t>0$, and moreover depends smoothly on $t$. In particular, there exists a smooth function $p$ depending on $t>0$ and $(x,\dot x)\in T^1B_0(1)$ such that
\[ \mathbb P\Big(H_t\in(T^1\phi)^{-1}(A)\Big)=\int_Ap_t(x,\dot x)\dif x\dif\dot x, \]
where the integral is considered with respect to Lebesgue measure (normalised however the reader pleases, up to introducing a constant in $p$). Since $p$ is smooth and $D_\varepsilon$ is compact, there exists a constant $\|p\|>0$ such that we have $p_t(x,\dot x)\leq\|p\|$ for all $t\in[a/2,2b]$, $x\in\phi(D_\varepsilon)$, $\dot x\in\mathbb S^{n-1}$. It means that for any such $t$ and $0<\delta\leq\varepsilon$,
\[ \mathbb P(H_t\in D_\delta) \leq \|p\|\cdot\operatorname{Vol}(\phi(D_\delta))\cdot\operatorname{Vol}(\mathbb S^{n-1}). \]
Let $\phi^*g$ be the metric on $B_0(1)$ induced by the identification with $\mathcal U$, seen as a $n\times n$ matrix. Using smoothness and compactness again, here exists a constant $C>0$ such that $\phi^*g$ is bounded above by $C\operatorname{id}$ in $D_\varepsilon$. In particular, given $0<\delta\leq\varepsilon$ and a point $y$ in $D_\delta$, there exists (by compactness) a point $x$ in $D$ and a smooth curve $\gamma$ of length at most $\delta$ with endpoints $x$ and $y$. The points of $\gamma$ belong to $D_\varepsilon$, so the Euclidean length of $\phi\circ\gamma$ is at most $C$ times the length of $\gamma$, which means that $\phi(y)$ is included in $\phi(D)+B_0(C\delta)$. All in all,
\[      \phi(D_\delta)
\subset \phi(D)+B_0(C\delta)
\subset \big(B_0(1)\cap\mathbb R^{n-2}\big)\times\big(B_0(C\delta)\cap\mathbb R^2\big) \]
for all $\delta$ small enough, and the Euclidean volume of $\phi(D_\delta)$ is bounded by $\delta^2$ up to a constant factor, which concludes.

Note that Proposition \ref{prop:codim2} is obviously optimal in terms of dimension. If any kinetic motion $(H,\dot H)$ were to avoid codimension 1 manifolds, then it wouldn't be able to reach the boundary of small balls around its initial point, so it would have to be stationary: $H_t=H_0$, $\dot H_t=0$. If one wants to show that a set that is not a submanifold is completely avoided by kinetic Brownian motion, the above shows that it would be enough for the $\delta$-fattenings of its bounded subsets to have volume of order $o(\delta)$.

\subsection{Proof of Lemmas \ref{lem:defULambda} and \ref{lem:MarkovLambdaA}}
\label{ssec:proofL12}

Suppose $(H,\dot H)$ is defined as in \eqref{eq:defH} and \eqref{eq:defdotH}. In particular, $\dot H$ is driven by a standard Brownian motion $W$ with values in $\mathcal H_d$. We want to show that the processes $U$, $\Lambda$ and $A$ are well-defined as long as $H$ has distinct eigenvalues, that they take values in the spaces $U_d(\mathbb C)$, $\mathcal H_d^\Delta$ and $\mathbb S(\mathcal H_d)$ respectively, and that $(\Lambda,A)$ satisfies the equations stated in Lemma \ref{lem:MarkovLambdaA}.

For the first few steps, we can actually work with a fixed realisation of $(H,\dot H)$. If $H_0$ has multiple eigenvalues then there is nothing to prove, so assume it is not the case. Writing
\[ \dot u_t = \dot u(U_t^*H_tU_t,U_t^*\dot H_tU_t) \]
and since $\dif U_t = U_t\dot u_t\dif t$ by definition, we can use the usual Picard-Lindelöf theorem to see that $U$ is uniquely well-defined for a maximal time interval $[0,T)$. Note that $\dot u(\Lambda,A)$ is in $\mathfrak u_d$ whenever it is well-defined, even if $\Lambda$ is not diagonal, which means that $U$ is in $U_d(\mathbb C)$ until $T$. This directly implies that $A$ takes values in the sphere $\mathbb S(\mathcal H_d)$, since $\dot H$ does and $K\mapsto U^*KU$ is an isometry of $\mathcal H_d$ for every unitary $U$.

If $U_t^*H_tU_t$, up to time $T$, stays uniformly away from the closed set of matrices with multiple diagonal entries, then $\dot u_t$ is bounded for $t\in[0,T)$, which means that $U_t$ converges to a limit as $t\to T$. Since $U_T^*H_TU_T$ has distinct diagonal entries, we can then apply Picard-Lindelöf again at $T$ and get a solution defined on a larger interval $[0,T+\varepsilon)$, which is a contradiction. Therefore, $T$ must be at least as large as the stopping time when at least two diagonal coefficients of $U_t^*H_tU_t$ become equal, or more precisely $T\geq\sup_{\varepsilon>0}T_\varepsilon$ for $T_\varepsilon$ the first time when $U_t^*H_tU_t$ is at distance at most $\varepsilon$ from the closed set of matrices with at least two diagonal coefficients being equal. Since $\dot u_t$ is not defined at this instant, we have in fact $T=\sup_{\varepsilon>0}T_\varepsilon$.

It implies that $\Lambda$ and $A$ are well-defined up to this time as well. Moreover, if we show that $\Lambda$ is diagonal up to $T$, then in fact the diagonal entries of $\Lambda_t=U_t^*H_tU_t$ are precisely the eigenvalues of $H_t$, and the collapse of the diagonal entries of the former corresponds to that of the eigenvalues of the latter, i.e. $T=\tau$.

Using the representation of $(H,\dot H)$ involving $W$, we see that for all $t<T$ we must have
\begin{align*}
\dif U_t & = U_t\dot u(\Lambda_t,A_t)\dif t, \\
\dif\Lambda_t & = \big(\dot u_t^*\Lambda_t + \Lambda_t\dot u_t\big)\dif t + A_t\dif t, \\
\dif A_t & = \big(\dot u_t^*A_t + A_t\dot u_t\big)\dif t
           + U_t^*\Big(\dif W_t - \dot H_t\Tr(\dot H_t^*\dif W_t) - \frac{d^2-1}2\dot H_t\dif t\Big)U_t \\
         & = \big(\dot u_t^*A_t + A_t\dot u_t\big)\dif t
           + U_t^*\dif W_tU_t - A_t\Tr(A_t^*U_t^*\dif W_tU_t) - \frac{d^2-1}2A_t\dif t.
\end{align*}
Since $U$ is unitary, the integral
\[ B:t\mapsto\int_0^tU_s^*\dif W_sU_s \]
defines a standard Brownian motion in $\mathcal H_d$, and $A$ satisfies the equation described in Lemma \ref{lem:MarkovLambdaA}. Moreover, if $\Lambda$ stays diagonal, then in fact
\[  \dif\Lambda_t
 = \pi^\Delta\big(\dif\Lambda_t\big)
 = \pi^\Delta\big(\dot u_t\Lambda_t + \Lambda_t\dot u_t\big)\dif t
 + \pi^\Delta(A_t)\dif t
 = \pi^\Delta(A_t)\dif t \]
according to equation \eqref{eq:productdotuLambda}, so $\Lambda$ satisfies the equation given in Lemma \ref{lem:MarkovLambdaA}. So the last thing we need to prove is that $\Lambda$ stays diagonal for all $t<T$.

This last fact is essentially a consequence of uniqueness for strong solutions of stochastic differential equations. Indeed, we can define $(\Lambda^B,A^B)$ as the solution of
\begin{align*}
\dif\Lambda^B_t & = \big((\dot u^B_t)^*\Lambda^B_t + \Lambda^B_t\dot u^B_t\big)\dif t + A^B_t\dif t, \\
\dif A^B_t & = \big((\dot u^B_t)^*A^B_t + A^B_t\dot u^B_t\big)\dif t
             + \dif B_t - A^B_t\Tr((A^B_t)^*\dif B_t) - \frac{d^2-1}2A^B_t\dif t
\end{align*}
for $\dot u^B = \dot u(\Lambda^B,A^B)$ and initial condition $(\Lambda^B,A^B)_0=(\Lambda,A)_0$, seen as a process with values in the open set of $\mathcal H_d^\Delta\times\mathcal H_d$ where the first component has distinct eigenvalues. It is defined on a (random) maximal interval $[0,T^B)$. The pairs $(\Lambda,A)$ and $(\Lambda^B,A^B)$ (or rather the inclusion of the latter in the full space $\mathcal H_d\times\mathcal H_d$) are solution to the same stochastic differential equation for all times before $T$ and $T^B$, so they are equal and $\Lambda$ is actually diagonal over $[0,T\wedge T^B)$. However, over the event $\{T^B<T\}$, the limit $(\Lambda^B,A^B)_{T^B}$ is well-defined in the large space where $(\Lambda,A)$ takes values, namely it is $(\Lambda,A)_{T^B}$ where $\Lambda_{T^B}\in\mathcal H_d$ with distinct diagonal entries and $A_{T^B}\in\mathcal H_d$. Since $\mathcal H_d^\Delta$ is closed, then in fact $(\Lambda^B,A^B)$ admits a limit in the small space as $t$ approaches $T^B$. But this event has measure zero according to the classical explosion criterion for equations with smooth coefficients, so the event $\{T^B<T\}$ has measure zero and $T\wedge T^B=T$, hence $\Lambda$ is diagonal for all times $t<T$. Though we don't need it, it also means that $(\Lambda,A)$ and $(\Lambda^B,A^B)$ are solutions to the same equation in the small space, and by maximality $T\leq T^B$, so we have in fact $T=T^B$ and $(\Lambda,A)=(\Lambda^B,A^B)$ always.

As discussed above, this concludes the proof of Lemmas \ref{lem:defULambda} and \ref{lem:MarkovLambdaA}.

\subsection{Projections of spherical Brownian motions}
\label{ssec:sphericalBM}

Let $X$ be a standard Brownian motion on the sphere $\mathbb S(\mathbb R^n)$, and $X^{[k]}$ its projection $(X^1,\ldots,X^k)$. The case we have in mind is $n=d^2$ and $k=d$. One way to define such an $X$ is to fix a standard Brownian motion $B$ with values in $\mathbb R^n$ and set $X$ the solution of
\[ \dif X_t
 = {}\circ\dif B_t - X_tX_t^*\circ\dif B_t
 = \dif B_t - X_tX_t^*\dif B_t - \frac{n-1}2X_t\dif t.  \]
We want to show that $X^{[k]}$ is Markovian. Since the projection $x\mapsto x^{[k]}$ is smooth, by Itô's formula we know that the process is solution to a stochastic differential equation of the form
\[  \dif X^{[k]}_t = b(X_t)\dif t + \sigma(X_t)\dif B_t  \]
with $b$ and $\sigma$ smooth. Such coefficients are uniquely determined, as they can be deduced from the generator $\frac12\Delta_\mathbb{S}$ of $X$ applied to functions of the form $x\mapsto f(x^{[k]})$. Moreover, for any fixed rotation $R$, since the law of $X$ is invariant under $R$, we know that
\[  \dif (RX)^{[k]}_t = (b\circ R)(X_t)\dif t + (\sigma\circ R)(X_t)\dif B_t. \]

For any fixed $x\in\mathbb S(\mathbb R^n)$, let $R$ be a rotation leaving the first $k$ coordinates invariant and such that $Rx = (x^1,\ldots,x^k,\sqrt{1-r^2},0,\ldots,0)$ with $r=\big|x^{[k]}\big|$. Then $(RX)^{[k]}=X^{[k]}$ for all times, so
\[ \dif X^{[k]}_t
 = b(X_t)\dif t + \sigma(X_t)\dif B_t
 = (b\circ R)(X_t)\dif t + (\sigma\circ R)(X_t)\dif B_t, \]
and by uniqueness
\[  b(x) = b(Rx) = b(x^1,\ldots,x^k,\sqrt{1-r^2},0,\ldots,0) = b^{[k]}\big(x^{[k]}\big)  \]
for a fixed function $b^{[k]}$ that is smooth away from $\{x^{[k]}=0\}$. We can use the same reasoning for $\sigma$, so that $X^{[k]}$ is Markovian, at least away from $\{x^{[k]}=0\}$.

Note that if $k\geq2$, which holds in our case, then almost surely, $X^{[k]}$ avoids zero except possibly for $t=0$. One way to see it is to notice that the traces of $X$ and $B/|B|$ have the same distribution, so that the probability that the trace of $X$ contains a point $x$ with $x^{[k]}=0$ is equal to that of $B$ containing a point $y$ with $y^{[k]}=0$. But $B^{[k]}$ is a standard Brownian motion in dimension at least 2, so it avoids zero almost surely.

Similar symmetry considerations can be used to show that $\theta:=X^{[k]}/|X^{[k]}|$ has to be a Brownian motion of the sphere, up to a time change depending only on the $r$ described above. In fact, using the decomposition of $\dif B_t$ along the following subspaces
\begin{align*}
E_t & := \{ x=(x^1,\ldots,x^k,0,\ldots,0)\text{ orthogonal to }X_t \} \\
F_t & := \{ x=(0,\ldots,0,x^{k+1},\ldots,x^n)\text{ orthogonal to }X_t \} \\
G_t & := (E_t\oplus F_t)^\perp = \mathbb R(X^1,\ldots,X^k,0\ldots,0) \oplus \mathbb R(0,\ldots,0,X^{k+1},\ldots,X^n),
\end{align*}
valid when $0\neq X^{[k]}\neq X$, one can show that
\begin{align*}
\dif(r^2)_t &= 2\sqrt{(1-r^2_t)r^2_t}\,\dif B^r_t + (k-nr^2_t)\dif t \\
\dif\theta_t &= \frac1{\sqrt{r_t^2}}(\dif B^\theta_t - \theta_t\theta_t^*\dif B^\theta_t) - \frac1{r_t^2}\frac{k-1}2\theta_t\dif t \\
\dif\phi_t &= \frac1{\sqrt{1-r_t^2}}(\dif B^\phi_t - \phi_t\phi_t^*\dif B^\phi_t) - \frac1{1-r_t^2}\frac{n-k-1}2\phi_t\dif t
\end{align*}
for $X=(r\theta,\sqrt{1-r^2}\phi)$ and $(B^r,B^\theta,B^\phi)$ a standard Brownian motion on $\mathbb R^{n+1}$. A (different) complete proof, as well as pathwise uniqueness, is described by Mijatović, Mramor and Uribe in \cite{MMU}.

\subsection{Proof of Lemma \ref{lem:MarkovPhi}}
\label{ssec:proofL3}

We want to show that $(\Lambda,\dot\Lambda)$ is Markovian if and only if the vector field $\Phi$ depends on $A$ only through its diagonal $\pi^\Delta(A)$. The indirect implication is clear: if $\Phi(\Lambda,A) = \overline\Phi(\Lambda,\pi^\Delta(A))$, then
\begin{align*}
\dif\Lambda_t &= \pi^\Delta(A_t)\dif t \\
\dif \pi^\Delta(A)_t
 &= \overline\Phi\big(\Lambda_t,\pi^\Delta(A)_t\big)\dif t
  + b^\Delta\big(\pi^\Delta(A_t)\big)\dif t + \sigma^\Delta\big(\pi^\Delta(A_t)\big)\dif B_t,
\end{align*}
where
\[ \dif X^\Delta_t = b^\Delta(X^\Delta_t)\dif t + \sigma^\Delta(X^\Delta_t)\dif B_t \]
is the equation describing the projection $X^\Delta = \pi^\Delta(X)$ on $\mathcal H_d^\Delta$ of a spherical Brownian motion $X$ in $\mathbb S(\mathcal H_d)$, as discussed in the previous section. Then $(\Lambda,\dot\Lambda)=(\Lambda,\pi^\Delta(A))$ is the solution of a self-contained SDE, so it is Markovian.

Conversely,  suppose that $(\Lambda,\dot\Lambda)$ is Markovian. Let $L^\Delta$ be its generator, $L$ that of $(\Lambda,A)$. For $f:\mathcal H_d^\Delta\times\mathcal H_d^\Delta\to\mathbb R$ regular enough, we should have
\[ L\big(f\circ(\operatorname{id},\pi^\Delta)\big)(\Lambda_0,A_0)
 = \frac{\dif}{\dif t}_{|t=0}\mathbb E_{\Lambda_0,A_0}[f(\Lambda_t,\pi(A_t))]
 = (L^\Delta f)\big(\Lambda_0,\pi^\Delta(A_0)\big). \]
For instance, one can see that this holds for $f$ smooth with compact support.

Let $X^\Delta=\pi^\Delta(X)$, as above, be the projection of a spherical Brownian motion on $\mathbb S(\mathcal H_d)$, and set $I$ its integral (i.e. $\dif I_t=X^\Delta_t\dif t$). According to the previous section, $(I,X^\Delta)$ is Markovian. In particular, $(I,X)$ and $(I,X^\Delta)$ both admit generators $\widetilde L$ and $\widetilde L^\Delta$, and they are linked by the same relation
\[ \widetilde L\big(f\circ(\operatorname{id},\pi^\Delta)\big)(I_0,X_0)
 = (\widetilde L^\Delta f)\big(I_0,\pi^\Delta(X_0)\big), \]
for instance when $f:\mathcal H_d^\Delta\times\mathcal H_d^\Delta\to\mathbb R$ is smooth with compact support.

As mentioned above, the only difference between $L$ and $\widetilde L$ is the additional vector field $\Phi$ acting on the second component:
\[ (L-\widetilde L)g(\Lambda_0,A_0)
 = D_{\!A}\,g(\Lambda_0,A_0)(\Phi(\Lambda_0,A_0))=:(\Phi\cdot\nabla_{\!A}\,g)(\Lambda_0,A_0), \]
where $D_{\!A}\,g$ is the differential of $g:\mathcal H_d^\Delta\times\mathcal H_d\to\mathbb R$ with respect to its second variable. In particular,
\[ \big(\Phi\cdot\nabla_A(f\circ(\operatorname{id},\pi^\Delta))\big)(\Lambda_0,A_0)
 = (L^\Delta f-\widetilde L^\Delta f)\big(\Lambda_0,\pi^\Delta(A_0)\big). \]
The right hand side depends on $A_0$ only through $\pi^\Delta(A_0)$, so the left hand side must be a function of $(\Lambda,\pi^\Delta(A))$. Since we can deduce a given vector field $\Psi$ with values in $\mathcal H_d^\Delta$ by the action of the operator $\Psi\cdot\nabla_{\!A}$ on functions of the form $f\circ(\operatorname{id},\pi^\Delta)$ (choose for instance a collection of $f$ smooth with compact support such that $f(\Lambda,\dot\Lambda) = \dot\Lambda_{ii}$ on a small open set), $\Phi$ actually factors through $(\Lambda,A)\mapsto(\Lambda,\pi^\Delta(A))$ as expected.

\subsection{The case \texorpdfstring{$d=2$}{d=2}}
\label{ssec:dim2}

We have seen at the end of Section \ref{ssec:criterion} that in dimension $d=2$, the process $(\Lambda,\dot\Lambda)$ is Markovian, using the fact that
\[ \Phi(\Lambda,A)_{11}
 = - \Phi(\Lambda,A)_{22}
 = \frac{|A_{12}|^2+|A_{21}|^2}{\Lambda_{11}-\Lambda_{22}}
 = \frac{1-|A_{11}|^2-|A_{22}|^2}{\Lambda_{11}-\Lambda_{22}}. \]
In fact, we can use this expression and the equation satisfied by $(\Lambda,A)$, given in Lemma \ref{lem:MarkovLambdaA}, to get the equation for the evolution. Write $\lambda$ and $\mu$ for the eigenvalues $\Lambda_{11}$ and $\Lambda_{22}$ of $H$. For $B$ a standard Brownian motion on $\mathcal H_2$, define the martingales
\[ M^\lambda : t \mapsto (B_{11})_t - \int_0^t\dot\lambda_s\Tr(A_s^*\dif B_s)
   \quad\text{ and }\quad
   M^\mu : t \mapsto (B_{22})_t - \int_0^t\dot\mu_s\Tr(A_s^*\dif B_s). \]
Then
\begin{align*}
\dif\lambda_t &= \dot\lambda_t\dif t, &
\dif\dot\lambda_t &= +\frac{1-\dot\lambda_t^2-\dot\mu_t^2}{\lambda_t-\mu_t}\dif t
                     + \dif M^\lambda_t
                     - \frac{d^2-1}2\dot\lambda_t\dif t, \\
\dif\mu_t &= \dot\mu_t\dif t, &
\dif\dot\mu_t &= -\frac{1-\dot\lambda_t^2-\dot\mu_t^2}{\lambda_t-\mu_t}\dif t
                     + \dif M^\mu_t
                     - \frac{d^2-1}2\dot\mu_t\dif t.
\end{align*}

Writing $A^\Re_{12}$ and $A^\Im_{12}$ for the real and imaginary parts of $A_{12}$, and similarly for the real and imaginary parts of $B_{12}$,
\begin{align*}
\dif M^\lambda_t
& = (1-\dot\lambda^2_t)\dif (B_{11})_t
  - \dot\lambda_t(\overline A_{12})_t\dif(B_{12})_t
  - \dot\lambda_t(\overline A_{21})_t\dif(B_{21})_t
  - \dot\lambda_t\dot\mu_t\dif(B_{22})_t \\
& = (1-\dot\lambda^2_t)\dif (B_{11})_t
  - 2\dot\lambda_t(A_{12}^\Re)_t\dif(B_{12}^\Re)_t
  - 2\dot\lambda_t(A_{12}^\Im)_t\dif(B_{12}^\Im)_t
  - \dot\lambda_t\dot\mu_t\dif(B_{22})_t, \\
\dif M^\mu_t
& = -\dot\lambda_t\dot\mu_t\dif (B_{11})_t
  - 2\dot\mu_t(A_{12}^\Re)_t\dif(B_{12}^\Re)_t
  - 2\dot\mu_t(A_{12}^\Im)_t\dif(B_{12}^\Im)_t
  + (1-\dot\mu_t^2)\dif(B_{22})_t.
\end{align*}
Since $\sum_{ij}|A_{ij}|^2=1$, we deduce $2|A_{12}^\Re|^2+2|A_{12}^\Im|^2 = 1-\dot\lambda^2-\dot\mu^2$, and we find the bracket of $M^\lambda$:
\[ \dif\langle M^\lambda,M^\lambda\rangle_t
 = (1-\dot\lambda_t^2)^2\dif t
 + \dot\lambda_t^2(1-\dot\lambda_t^2 - \dot\mu_t^2)\dif t
 + \dot\lambda_t^2\dot\mu_t^2\dif t
 = (1 - \dot\lambda_t^2)\dif t. \]
Similarly the bracket of $M^\mu$ grows as $(1-\dot\mu_t^2)\dif t$. The covariance term is given by
\[ \dif\langle M^\lambda,M^\mu\rangle_t
 = - \dot\lambda_t\dot\mu_t(1-\dot\lambda_t^2)\dif t
   + \dot\lambda_t\dot\mu_t(1-\dot\lambda_t^2-\dot\mu_t^2)\dif t
   - \dot\lambda_t\dot\mu_t(1-\dot\mu_t^2)\dif t
 = - \dot\lambda_t\dot\mu_t\dif t, \]
as stated in Lemma \ref{lem:casesond}.

Note that it corresponds to the diffusion term for the projection of a spherical Brownian motion, as described in \cite{MMU}. Indeed, as explained in the proof outline, the only difference between $A$ and a spherical Brownian motion is a drift term. Alternatively, one can also study the trace and determinant of $H$, and deduce the process satisfied by the eigenvalues, since they are the roots of the polynomial $X^2-\Tr(H)X+\operatorname{det}(H)$.

\subsection{The case \texorpdfstring{$d\geq3$}{d≥3}}
\label{ssec:dim3+}


We show that in this case, the vector field $\Phi(\Lambda,A)$ depends on the off-diagonal elements of $A$. As stated in Lemma \ref{lem:MarkovPhi}, this will show that $(\Lambda,\dot\Lambda)$ cannot be Markovian. We will use the expression given in equation \eqref{eq:concisePhi}.

In dimension 3, it is a direct computation to see that for any $\Lambda$ with distinct eigenvalues, the following two choices for $A$ give different $\Phi(\Lambda,A)_{11}$, although they are equal on the diagonal:
\begin{align*}
A &= \begin{pmatrix} 0 & 1 & 0 \\ 1 & 0 & 0 \\ 0 & 0 & 0 \end{pmatrix}, &
\widetilde A &= \begin{pmatrix} 0 & 0 & 0 \\ 0 & 0 & 1 \\ 0 & 1 & 0 \end{pmatrix}.
\end{align*}
In fact, one gets $\Phi(\Lambda,A)_{11} = 2/(\Lambda_{11}-\Lambda_{22})$, whereas $\Phi(\Lambda,\widetilde A)_{11} = 2/(\Lambda_{11}-\Lambda_{33})$. In higher dimension, chose $A$ and $\widetilde A$ to be zero except on the top left $3\times3$ minor, which is given by the above expressions.

One might be interested in weaker versions of Markovian behaviour. For instance, one may try to restrict to a well-chosen family of initial conditions and hope that $(\Lambda,\dot\Lambda)$ be Markovian up to some exit time. Let $N$ be some submanifold of $\mathcal H_d^\Delta\times\mathbb S(\mathcal H_d)$, and consider the process $(\Lambda,A)^\dagger$ as $(\Lambda,A)$ stopped when it exits $N$, with values in the one-point compactification $N\cup\{\dagger\}$ of $N$. By convention, $\dagger$ is absorbing. We will show that this localisation of $(\Lambda,A)$ does not induce a Markovian $(\Lambda,\dot\Lambda)$, except in the uninteresting case $(\Lambda,A)^\dagger_t=\dagger$ for all $t>0$, where the motion was stopped immediately.

If $N$ is not of maximal dimension, then by hypoellipticity of $(\Lambda,A)$ (see the governing equation in Lemma \ref{lem:MarkovLambdaA}), the unstopped $(\Lambda,A)_t$ admits a smooth density for arbitrarily small $t>0$, and the probability that $(\Lambda,A)_t$ belongs to $N$ is zero ($N$ has zero Lebesgue measure as a submanifold of positive codimension). In particular $(\Lambda,A)^\dagger_t=\dagger$ for all $t>0$.

We consider now the case where $N$ is of maximal dimension; in other words, $N$ is an open set. Unwinding the proof of Lemma \ref{lem:MarkovPhi} in Section \ref{ssec:proofL3}, we see that the induced $(\Lambda,\dot\Lambda)$ will be Markovian if and only if the vector field $\Phi(\Lambda,A)$ factors through $(\Lambda,A)\mapsto(\Lambda,\pi^\Delta(A))$ when restricted to $N$.

Since $N$ is open, we can choose $(\Lambda,A)\in N$ such that $\Lambda$ has distinct eigenvalues and $A$ is not diagonal. There exists $i\neq j$ such that $A_{ij}\neq0$; moreover, because $d\geq3$, there exists $k\notin\{i,j\}$. Choose a unit complex number $u$ such that $A_{ik}=|A_{ik}|u$. Set $A^\varepsilon\in\mathcal H_d$ equal to $A$ everywhere but at $(i,j)$, $(i,k)$ and their symmetric counterparts $(j,i)$ and $(k,i)$, which are instead defined by
\begin{align*}
A^\varepsilon_{ij} & = A_{ij}\cos(\varepsilon), &
A^\varepsilon_{ik} & = A_{ik} + \mathsf i|A_{ij}|\sin(\varepsilon)u_{ik}. &
\end{align*}
The important features of this perturbation are that $|A^\varepsilon_{ij}|$ is not constant around $\varepsilon=0$, and $A^\varepsilon$ stays on the sphere $\mathbb S(\mathcal H_d)$; indeed,
\[ |A^\varepsilon_{ij}|^2 + |A^\varepsilon_{ik}|^2
 = |A_{ij}|^2\cos(\varepsilon)^2 + |A_{ik}|^2 + |A_{ij}|^2\sin(\varepsilon)^2
 = |A_{ij}|^2 + |A_{ik}|^2. \]
In particular, the entries of $\Phi(\Lambda,A^\varepsilon)$ along the diagonal should not depend on $\varepsilon$. However, writing $\lambda^\ell:=\Lambda_{\ell\ell}$ for the $\ell$th eigenvalue of $\Lambda$,
\begin{align*}
\frac{\dif^2}{\dif\varepsilon^2}\big(\Phi(\Lambda,A^\varepsilon)_{ii}\big)
& = \frac2{\lambda^i-\lambda^j}\cdot\frac{\dif^2}{\dif\varepsilon^2}|A^\varepsilon_{ij}|^2
  + \frac2{\lambda^i-\lambda^k}\cdot\frac{\dif^2}{\dif\varepsilon^2}|A^\varepsilon_{ik}|^2 \\
& = \left(\frac2{\lambda^i-\lambda^j} - \frac2{\lambda^i-\lambda^k}\right)
      \cdot\frac{\dif^2}{\dif\varepsilon^2}|A^\varepsilon_{ij}|^2 \\
& = \frac{\lambda^k-\lambda^j}{(\lambda^i-\lambda^j)(\lambda^i-\lambda^k)}\cdot4|A^\varepsilon_{ij}|^2\cos(2\varepsilon).
\end{align*}
The first factor is well-defined and non zero since $\Lambda$ has distinct eigenvalues, and the second is non-zero for $\varepsilon=0$. This shows that $\Phi(\Lambda,A)_{ii}$ depends on the off-diagonal entries of $A$ over $N$, hence $(\Lambda,\dot\Lambda)$ cannot be Markovian according to Lemma \ref{lem:MarkovPhi}.

\end{document}